\documentclass[12pt]{article}

\usepackage{graphicx}

\def\ra{\rightarrow}

\def\ss{\subseteq}

\def\Re{\hbox{\rm Re}\,}
\def\Im{\hbox{\rm Im}\,}

 \def\HollowBox #1#2{{\dimen0=#1 \advance\dimen0 by -#2       
       \dimen1=#1 \advance\dimen1 by #2                       
        \vrule height #1 depth #2 width #2                    
        \vrule height 0pt depth #2 width #1                   
        \llap{\vrule height #1 depth -\dimen0 width \dimen1}%
       \hskip -#2                                             
       \vrule height #1 depth #2 width #2}}                   
 \def\BoxOpTwo{\mathord{\HollowBox{6pt}{.4pt}}\;}             

\def\endpf{\hfill $\BoxOpTwo$}

\font\teneufm=eufm10
\font\seveneufm=eufm7
\font\fiveeufm=eufm5
\newfam\eufmfam
\textfont\eufmfam=\teneufm
\scriptfont\eufmfam=\seveneufm
\scriptscriptfont\eufmfam=\fiveeufm

\newfam\msbfam
\font\tenmsb=msbm10  scaled \magstep1 \textfont\msbfam=\tenmsb
\font\sevenmsb=msbm7 scaled \magstep1 \scriptfont\msbfam=\sevenmsb
\font\fivemsb=msbm5  scaled \magstep1 \scriptscriptfont\msbfam=\fivemsb
\def\Bbb{\fam\msbfam \tenmsb}

\def\CC{{\Bbb C}}

\def\NN{{\Bbb N}}

\newtheorem{theorem}{Theorem}

\newtheorem{proposition}[theorem]{Proposition}

\newtheorem{definition}{Definition}
\newtheorem{example}[definition]{EXAMPLE}

\makeindex

\begin{document}

\begin{center}
\Large \bf On Limits of Sequences of Holomorphic Functions
\end{center}

\begin{center}
\large Steven G. Krantz\footnote{Author supported in part by a grant from the
National Science Foundation and a grant from the Dean of the Graduate School
at Washington University in St.\ Louis.}\footnote{{\bf Key Words:}  holomorphic 
function, pointwise limit, uniform limit, elliptic equations.}\footnote{{\bf MR Subject Classification Numbers:}
30B60, 30E10, 30E99, 32A30, 32A99.}
\medskip \\
\small \today

\end{center}
\vspace*{.18in}

\begin{quote} \sl
{\bf Abstract:}  We study functions which are the pointwise limit of a sequence
of holomorphic functions.  In one complex variable this is a classical topic,
though we offer some new points of view and new results.  Some novel
results for solutions of elliptic equations will be treated.  In several complex
variables the question seems to be new, and we explore some new avenues.
\end{quote}

\setcounter{section}{-1}

\section{Introduction}

It is a standard and well known fact from complex function theory (which
appears to be due to Stieltjes [STE], although see also Vitali's theorem
in [TIT] and Weierstrass's complete works [WEI]) that if $\{f_j\}$ is a sequence of
holomorphic functions on a planar domain $\Omega$ and if the sequence
converges {\it uniformly on compact subsets of $\Omega$} then the limit
function is holomorphic on $\Omega$. Certainly this result is one of
several justifications for equipping the space of holomorphic functions on
$\Omega$ with the compact-open topology (see also [LUR], where this point
of view is developed in detail from the perspective of functional
analysis).

Considerably less well known is the following result of William Fogg Osgood [OSG]:

\begin{theorem} \sl
Let $\{f_j\}$ be a sequence of holomorphic functions on a planar
domain $\Omega$.  Assume that the $f_j$ converge {\it pointwise} to
a limit function $f$ on $\Omega$.  Then $f$ is holomorphic on a dense,
open subset of $\Omega$.   The convergence is uniform on compact
subsets of the dense, open set.
\end{theorem}

This result is not entirely obvious; it is certainly surprising and
interesting.  For completeness, we now offer a proof of the theorem:
\medskip \\

\noindent {\bf Proof of the Theorem:}  Let $U$ be a nonempty open subset
of $\Omega$ with compact closure in $\Omega$.  Define, for 
$k = 1, 2, \dots$,
$$
S_k = \{z \in \overline{U}: |f_j(z)| \leq k \ \hbox{for all} \ j \in \NN\} \, .
$$
Since the $f_j$ converge at each $z \in \overline{U}$, certainly the set
$\{f_j(z): j \in \NN\}$ is bounded for each fixed $z$.  
So each $z \in \overline{U}$ lies in some $S_k$.
In other words,
$$
\overline{U} = \bigcup_k S_k \, .
$$

Now of course $\overline{U}$ is a complete metric space (in the
ordinary Euclidean metric), so the Baire category theorem tells
us that some $S_k$ must be ``somewhere dense'' in $\overline{U}$.
This means that $\overline{S_k}$ will contain a nontrivial Euclidean metric
ball (or disc) in $\overline{U}$.  Call the ball ${\cal B}$.  Now it is a simple
matter to apply Montel's theorem on ${\cal B}$ to find a subsequence
$f_{j_k}$ that converges uniformly on compact sets to a limit
function $g$.  But of course $g$ must coincide with $f$, and
$g$ (hence $f$) must be holomorphic on ${\cal B}$.  

Since the choice of $U$ in the above arguments was arbitrary, the conclusion
of the theorem follows.
\endpf
\smallskip \\

\noindent {\bf Remark:} An alternative approach, which avoids the explicit
use of Montel's theorem, is as follows. Once one has identified an $S_k$
whose closure contains a ball or disc $D(P,r)$, let $\gamma:[0,1] \rightarrow \Omega$ be a simple,
closed, rectifiable curve in $D(P,r)$. Then of course the image $\widetilde{\gamma} \equiv \{\gamma(t): t \in [0,1]\}$ of
$\gamma$ is a compact set. Let $\epsilon > 0$. By Lusin's theorem, the
sequence $f_j$ converges uniformly on some subset $E \ss \widetilde{\gamma}$ with
the property that the linear measure of $\widetilde{\gamma} \setminus E$ is less
than $\epsilon$.  Let $K$ be a compact subset of the open region surrounded by $\gamma$,
and let $\delta > 0$ be the Euclidean distance of $K$ to $\widetilde{\gamma}$. 
Let $\epsilon^* > 0$ and choose $J > 0$ so large that when $\ell, m > J$ then
$$
|f_\ell(z) - f_m(z)| < \epsilon^*
$$
for all $z \in E$.  Then, for $w \in K$,
\begin{eqnarray*}
|f_\ell(w) - f_m(w)| & = & \left | \frac{1}{2\pi i} \oint_\gamma \frac{f_\ell(\zeta) - f_m(\zeta)}{\zeta - w} \, d\zeta \right | \\
                 & \leq & \frac{1}{2\pi} \int_E \frac{\epsilon^*}{\delta} \, ds +
		           \frac{1}{2\pi} \int_{\widetilde{\gamma} \setminus E} \frac{2k}{\delta} \, ds \\
		 & \leq & \frac{\epsilon^*\cdot |E|}{2\pi \cdot \delta} + \frac{\epsilon \cdot 2k}{2\pi \cdot \delta} \, . \\
\end{eqnarray*}
Thus we see that we have uniform convergence on $K$.  And the holomorphicity follows as usual.
\endpf 
\smallskip \\

The book [REM] contains a nice treatment of some of the one-complex-variable theory
related to Osgood's theorem.

The next example is inspired by ideas in [ZAL, pp.\ 131--133].  It demonstrates
that Osgood's theorem has substance, and describes a situation that
actually occurs.  A thorough discussion of many of the ideas treated here---from a somewhat
different point of view---appears in [BEM].  In fact [BEM] presents quite a different
contruction of an example that illustrates Theorem 1.

\begin{example} \rm
Let 
$$
U = \{z \in \CC: |\Re z| < 1, |\Im z| < 1\} \, .
$$
For $j = 1, 2, \dots$, define
$$
S_j = \{z \in U: \Re z = 0 \ \ \hbox{or} \ \ \Im z = 0, |\Re z| \leq 1 - 1/[j+2] \ , \
|\Im z| \leq 1 - 1/[j+2]\} \, .
$$
Also define
$$
T_j = \{z \in U: 1/[j+2] \leq |\Re z| \leq 1 - 1/[j+2], 1/[j+2] \leq |\Im z| \leq 1 - 1/[j+2]\} \, .
$$
We invite the reader to examine Figure 1 to appreciate these sets.

\begin{figure}
\centering				     
\includegraphics[height=3in, width=2.75in]{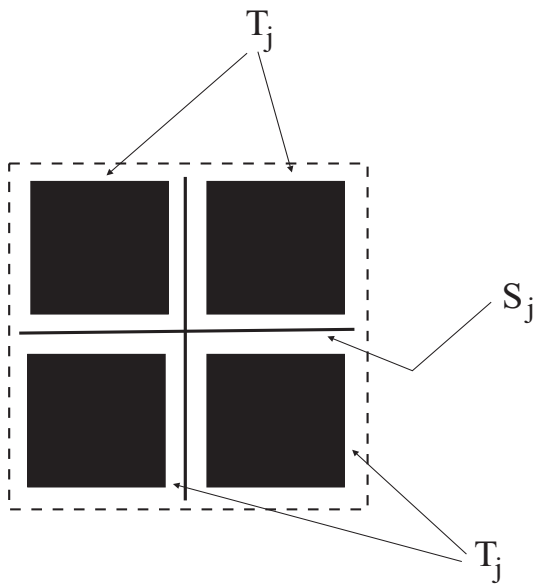}
\caption{The sets $S_j$ and $T_j$.}
\end{figure}

Now, for each $j$, we apply Runge's theoreom on $S_j \cup T_j$.  Notice that the complement of $S_j \cup T_j$ 
is connected, so that we can push the poles of the approximating functions to the complement
of $U$.  We are able then to produce for each $j$ a holomorphic function $f_j$ on $U$ such that
$$
|f_j(z) - 0| < \frac{1}{j} \quad \hbox{for} \ \ z \in T_j \, ,
$$
$$
|f_j(z) - 1| < \frac{1}{j} \quad \hbox{for} \ \ z \in S_j \, .
$$
Then it is easy to see that the sequence $\{f_j\}$ converges pointwise to the function
$f$ given by
$$
f(z) = \left \{ \begin{array}{lcr}
            0 \ & \ \hbox{if} \ & \ z \in U \setminus \{z \in U: \Re z = 0 \ \hbox{or} \ \Im z = 0\} \\
	    1 \ & \ \hbox{if} \ & \ z \in U \cap \{z \in U: \Re z = 0 \ \hbox{or} \ \Im z = 0\} \, .
	    	\end{array}
       \right.
$$
Thus the limit function $f$ is holomorphic on a dense open subset of $U$, and the exceptional set
is the two axes in $U$.
\endpf
\end{example}

One might ask what more can be said about the open, dense set $V$ on
which the limit function $f$ is holomorphic.  Put in other words,
what can one say about $\Omega \setminus V$?  Lavrentiev [LAV] was
the first to give a characterization of those open sets on which a pointwise
convergent sequence of holomorphic functions can converge.  Siciak [SIC]
has given a rather different answer in the language of capacity theory.  

In fact a suitable version of the theorem is true for harmonic functions, or more
generally for solutions of a uniformly elliptic partial differential equation of
second order.  We shall prove such results later in the present paper.	See also
[DAVI1], [NIC1], [NIC2], [STE] for related results.

It is a pleasure to thank T. W. Gamelin, D. Minda, D. Sarason,
and L. Zalcman for helpful discussions of the topics of this
paper. 
\bigskip \\

\section{More Results in the Classical Setting}

Our first new result for planar domains concerns harmonic functions:

\begin{theorem} \sl
Let $\{f_j\}$ be a sequence of harmonic functions on a planar
domain $\Omega$.  Assume that the $f_j$ converge {\it pointwise} to
a limit function $f$ on $\Omega$.  Then $f$ is harmonic on a dense
open subset of $\Omega$.
\end{theorem}

\noindent {\bf Sketch of the Proof of the Theorem:}  Proceed as
in the proof of the result for holomorphic functions.  It is
certainly true that a collection of harmonic functions on a planar
domain that is uniformly bounded on compacta will have a subsequence that converges
uniformly on compact sets.   This follows from easy estimates on the Poisson
kernel.  The rest of the argument is the same as before.
\endpf
\smallskip \\

\noindent {\bf Remark:}  A result of this nature is already contained in [GAG, Theorem 1].

\begin{theorem} \sl
Let ${\cal L}$ be a uniformly elliptic
operator, having locally bounded, $C^2$ coefficients, of order 2 on a planar domain $\Omega$. Let $\{f_j\}$
be a sequence of functions that are annihilated by ${\cal L}$
on $\Omega$. Assume that the $f_j$ converge {\it pointwise} to
a limit function $f$ on $\Omega$. Then $f$ is annihilated by ${\cal L}$ on a
dense open subset of $\Omega$. 
\end{theorem}
{\bf Proof:}  The proof is the same as the last result.  The only
thing to check is that a collection of functions annihilated by ${\cal L}$
that is bounded on compact sets will have a subsequence that converges uniformly
on compact sets.  This will follow, as in the harmonic case, from the Poisson
formula for ${\cal L}$ (see [GAR]).  The rest of the argument is the same.
\endpf
\smallskip \\

\begin{theorem} \sl
Let $\{f_j\}$ be a sequence of holomorphic functions on a planar
domain $\Omega$.  Suppose that there is a constant $M > 0$ such
that $|f_j(z)| \leq M$ for all $j$ and for all $z \in \Omega$.  Assume that the $f_j$ converge {\it pointwise} to
a limit function $f$ on $\Omega$.  Then $f$ is holomorphic on all
of $\Omega$.
\end{theorem}
\vspace*{.12in}

\noindent {\bf Remark:}  Of course the new feature in this last theorem
of Stieltjes [STE] is that we are assuming that the family $\{f_j\}$ is uniformly bounded.
This Tauberian hypothesis gives a stronger conclusion.  The
proof will now be a bit different.
\medskip \\

\noindent {\bf Proof:}  Let $U$ be an open subset of $\Omega$.
Then the argument from the proof of Theorem 1 applies immediately
on $U$.  Thus the limit function is holomorphic on $U$.  Since the
choice of $U$ was arbitrary, we are finished.
\endpf 
\smallskip \\

In fact there is a much weaker condition (than in the last theorem) that
will give the same result:

\begin{theorem}  \sl
Let $\{f_j\}$ be a sequence of holomorphic functions on a planar
domain $\Omega$.  Suppose that there is a nonnegative, integrable function
$g$ on $\Omega$ such that 
$|f_j (z)| \leq g(z)$ for all $j$ and for all $z \in \Omega$.  Assume that the $f_j$ converge {\it pointwise} to
a limit function $f$ on $\Omega$.  Then $f$ is holomorphic on all
of $\Omega$.
\end{theorem}
{\bf Proof:}  The proof is simplicity itself. Suppose that $K$ is a compact subset of $\Omega$.  Let $\varphi$ be a $C_c^\infty$
function on $\Omega$ that is identically equal to 1 on $K$.  Then
$$
f_j(z) = \frac{1}{\pi} \int \!\!\! \int \frac{f_j(\zeta) \overline{\partial} \varphi(\zeta)}{z - \zeta} \, dA(\zeta) 
$$
for each $j$.  Here $dA$ is Lebesgue area measure on $\CC$.  As a consequence,
$$
|f_j(z)| \leq \frac{1}{\pi} \int \!\!\! \int \frac{g(\zeta) |\overline{\partial} \varphi(\zeta)|}{|z - \zeta|} \, dA(\zeta) \, .
$$
Since $K$ has positive distance from the support of $\overline{\partial} \varphi$, we may now conclude
that the $\{f_j\}$ are uniformly bounded on $K$.  Theorem 4 therefore tells us that $f$ is holomorphic
on all of $\Omega$.
\endpf
\smallskip \\

The following result is discussed but not proved in [DAV]:

\begin{theorem} \sl
Let $\{f_j\}$ be a sequence of Schlicht functions on the unit disc
$D$ that converges pointwise.  Then the limit function is holomorphic
on all of $D$.
\end{theorem}
{\bf Proof:}  By the ``growth theorem'' (see [DUR]), any Schlicht function $f$ satisfies
$$
|f(z)| \leq |z| \cdot (1 - |z|)^{-2} 
$$
for all $z$ in the disc $D$.  It follows that the $f_j$ are uniformly bounded on compact
subsets of $D$.  So they form a normal family.  Thus there is a subsequence $\{f_{j_k}\}$
that converges uniformly on compact sets to some limit function $g$.  That function
$g$ is of course holomorphic everywhere.  But it must coincide with the pointwise
limit function.
\endpf 
\smallskip \\
	   
It is easy to see that results of the kind we are discussing here cannot
be true in the category of real analytic functions.  Indeed, the
Weierstrass approximation theorem tells us that {\it any} continuous
function is the uniform limit on compact sets of polynomials, hence of 
real analytic functions.  But we do have the following modified result:

\begin{proposition} \sl
Let $\{f_j\}$ be a sequence of real analytic functions on the
bounded interval $(a,b)$.  For each $\ell = 0, 1, 2, \dots$, let
$g^{(\ell)}$ denote the $\ell^{\rm th}$ derivative of the function $g$.
Assume that, for each $\ell$, there are positive constants $K$ and $R$ so
that, for all $j$ and $\ell$,
$$
|f_j^{(\ell)}| \leq K \cdot \frac{\ell !}{R^\ell} \, .
$$
Further assume that the sequence $\{f_j\}$ converges pointwise to
some function $f$ on $(a,b)$.  Then $f$ is real analytic on $(a,b)$.
\end{proposition}

\noindent The proposition bears out the heuristic that a sequence of
real analytic functions converging to a non-analytic function
like $f(x) = |x|$ on the interval $(-1,1)$ must have derivatives blowing
up.  Note that, for simplicity, we have stated the result in one dimension.
But it clearly holds in any dimension.
\smallskip \\

\noindent {\bf Proof of the Proposition:}  Without loss of generality, assume
that $a = -b$, so that the point 0 lies in the center of the domain interval.
Call the interval now $(-a,a)$.

Fix a positive integer $N$.  With a simple diagonalization argument, we may
choose a subsequence $f_{j_k}$ so that $f_{j_k}^{(\ell)}(0)$ converges
to some number $\alpha_\ell$ as $k \ra \infty$, each $\ell = 0, 1, 2, \dots, N$.  For each $k$, we may write
$$
f_{j_k}(x) = \sum_{\ell = 0}^N \frac{f_{j_k}^{(\ell)}}{\ell!} \cdot x^\ell + {\cal O}(x^{N+1}) \, .
$$
Letting $k \ra +\infty$, we find that
$$
f(x) = \sum_{\ell = 0}^N \frac{\alpha_\ell}{\ell!} \cdot x^\ell + {\cal O}(x^{N+1}) \, .
$$
Of course a similar identity holds at each point of $(-a,a)$ (not
just at 0).  It follows then, from the converse to the Taylor theorem (see [ABR],
[KRA2]), that $f$ is $C^\infty$ and $f^{(\ell)}$ at each point $x$
is given as the limit of a subsequence of the $f_j^{(\ell)}$ at that point.
Since the argument applies to show that every subsequence of $\{f_j\}$ has a subsequence
with this property, we may conclude that
$$
\lim_{j \ra \infty} f_j^{(\ell)}(x) = f^{(\ell)}(x)
$$
for each $x$ in the $(-a,a)$ and each $\ell = 0, 1, 2, \dots$.

Now we have assumed that
$$
|f_j^{(\ell)}(x)| \leq K \cdot \frac{\ell!}{R^\ell}
$$
for $x$ in a compact subset of $(-a,a)$.  Letting $j \ra \infty$ gives
$$
|f^{(\ell)}(x)| \leq K \cdot \frac{\ell!}{R^\ell} \, .
$$
We conclude then that $f$ is real analytic, as desired.
\endpf
\smallskip \\

\section{Results in Several Complex Variables}

The first result in $\CC^n$ is as follows.

\begin{theorem}  \sl
Let $\{f_j\}$ be a sequence of holomorphic functions on a 
domain $\Omega \ss \CC^n$.  Assume that the $f_j$ converge {\it pointwise} to
a limit function $f$ on $\Omega$.  Then $f$ is holomorphic on a dense
open subset of $\Omega$.  Also the convergence is uniform on compact
subsets of the dense open set.
\end{theorem}
{\bf Proof:}  The argument is the same as that for Theorem 1.
We need only note that Montel's theorem is still valid.  The
rest of the argument is the same.
\endpf
\smallskip \\

\noindent {\bf Remark:}  Just as in the Remark following the proof of Theorem
1, we could use the Henkin-Ramirez integral formula on small balls (see [KRA1, Ch.\ 8])
to give an alternative proof of this result.
\endpf
\smallskip \\

\begin{theorem} \sl
Let $\{f_j\}$ be a sequence of holomorphic functions on a 
domain $\Omega \ss \CC^n$.  Assume that the $f_j$ converge {\it pointwise} to
a limit function $f$ on $\Omega$.   Let $\ell$ be any complex
line in $\CC^n$.   Then the limit function $f$ is holomorphic on
a dense open subset of $\ell \cap \Omega$.
\end{theorem}
{\bf Proof:}  Of course we simply apply the argument from
the proof of Theorem 1 on $\ell \cap \Omega$.
\endpf
\smallskip \\

\noindent {\bf Remark:}  This is a stronger result than Theorem 8.  One may
note that something similar could be proved with ```complex line'' replaced
by ``complex analytic variety''.  It is not clear what the optimal result might be.

We note that [DAVI1] contains some results which characterize those sets whose characteristic function
is the pointwise limit of a sequence of holomorphic functions.  The main result of [wDAVI1] was
anticipated in the paper [DAM].  See also [DAVI2].  
\endpf
\smallskip \\

With a little effort, one may produce results in several complex variables that
introduce a new way to think about these theorems.  An example is this:

\begin{theorem} \sl
Let $\{f_j\}$ be a sequence of holomorphic functions on a domain
$\Omega \ss \CC^n$, $n \geq 2$.  Suppose that there is a function $f$ on $\Omega$ such that, for each analytic
disc $\varphi: D \ra \Omega$, the sequence $\{f_j \circ \varphi\}$ converges
uniformly on $\overline{D}$ to $f \circ \varphi$.  
Then $f$ is holomorphic on $\Omega$.
\end{theorem}
{\bf Proof:}   Fix $z \in \Omega$ and $c > 0$ small.  Fix an index $j$ and let $\varphi (\zeta) = (0,0, \dots, c \cdot \zeta, 0, \dots 0)$, where
the $\zeta$ appears in the $j^{\rm th}$ position.  Then our hypothesis says that
$f_j \circ \varphi (\zeta) = f_j((0,0, \dots, \zeta, 0, \dots 0)$ converges uniformly, for $\zeta \in D$ to $f$.
But the limit function will of course be holomorphic.  So $f$ is holomorphic in each variable separately.  It follows
then from Hartogs's theorem (see [KRA1]) that $f$ is genuinely holomorphic as a function of $n$ variables.
\endpf 
\smallskip \\

\section{Concluding Remarks}

It is clear that there is more to learn in the several complex variable setting.
We would like a result that has a chance of being sharp, so that the exceptional
set for convergence can be characterized (as in [SIC] for one complex variable).
This matter will be explored in future papers.
\bigskip \bigskip \bigskip \\

It is a pleasure to thank the referee for many helpful comments and suggestions.

\vfill
\eject

\noindent {\Large \sc References}
\smallskip \\

\begin{enumerate}

\item[{\bf [ABR]}]  R. Abraham and J. Robbin, {\it Transversal Mappings and Flows}, Benjamin,
New York, 1967.

\item[{\bf [BEM]}] A. Beardon and D. Minda, On the pointwise
limit of complex analytic functions, {\it Am.\ Math.\ Monthly}
110(2003), 289--297.

\item[{\bf [DAM]}] A. A. Danielyan and S. N. Mergelyan,
Sequences of polynomials converging on sets of type
$F_\sigma$, {\it Akad.\ Nauk.\ Armyan}, SSR Dokl.\ 86(1988),
54--57.

\item[{\bf [DAV]}] K. R. Davidson, Pointwise limits of analytic
functions, {\it Am.\ Math.\ Monthly} 90(1984), 391--394.
				     
\item[{\bf [DAVI1]}]  A. M. Davie, Pointwise limits of analytic
functions, {\it J. London Math.\ Soc.} 75(2007), 133--145.

\item[{\bf [DAVI2]}]  A. M. Davie, addendum to Pointwise limits of analytic
functions, {\it Proc.\ London Math.\ Soc.} 79(2009), 272.

\item[{\bf [DUR]}] P. Duren, {\it Univalent Functions},
Springer-Verlag, New York, 1983.

\item[{\bf [GAR]}]  P. Garabedian, {\it Partial Differential Equations}, $2^{\rm nd}$ ed., 
Chelsea Publishing, New York, 1986.

\item[{\bf [GAG]}] S. Gardiner and A. Gustafsson, Pointwise
convergence and radial limits of harmonic functions, {\it
Israeli J. Math.} 145(2005), 243--256.

\item[{\bf [GRK]}] R. E. Greene and S. G. Krantz, {\it Function
Theory of One Complex Variable}, $3^{\rm rd}$ ed., American
Mathematical Society, Providence, RI, 2006.

\item[{\bf [KRA1]}] S. G. Krantz, {\it Function Theory of
Several Complex Variables}, $2^{\rm nd}$ ed., American
Mathematical Society, Providence, RI, 2002.

\item[{\bf [KRA2]}] S. G. Krantz, Lipschitz spaces, smoothness of
functions, and approximation theory, {\it Expositiones Math.} 3(1983),
193-260.

\item[{\bf [KRP]}]  S. G. Krantz and H. R. Parks, {\it A Primer
of Real Analytic Functions}, $2^{\rm nd}$ ed., Birkh\"{a}user Publishing, Boston,
2002.

\item[{\bf [LAV]}] M. Lavrentiev, Sur les fonctions d'une
variable complexe repr\'{e}sentables par des s\'{e}ries de
polynomes, {\it Actualit\'{e}s Scientifiques et Industrielles}
441(1936), Hermannn \& Cie, Paris.

\item[{\bf [LUR]}] D. Luecking and L. Rubel, {\it Complex
Analysis: A Functional Analysis Approach}, Springer-Verlag,
New York, 1984.

\item[{\bf [MON]}] P. Montel, {\it Le\c{c}ons sur les Families
Normales de fonctions analytiques et leur applications},
Gauthier-Villars, Paris, 1927.

\item[{\bf [NIC1]}] C. P. Niculescu, Applications of elliptic
operator methods to $C^\infty$ convergence problem, {\it Rev.\
Roum.\ Math.\ Pures Appl.} 44(1999), 793--798.

\item[{\bf [NIC2]}] C. P. Niculescu, Function spaces attached
to elliptic operators, {\it Proc.\ of the National Conference
on Mathematical Analysis and Applications}, Timosoara,
December 12--13, 2000, pp.\ 239--250, University of West
Timosoara, 2000.		   

\item[{\bf [OSG]}] W. F. Osgood, Note on the functions defined
by infinite series whose terms are analytic functions of a
complex variable, with corresponding results for definite
integrals, {\it Ann.\ Math.} 3(1901), 25--34.

\item[{\bf [REM]}]  R. Remmert, {\it Classical Topics in Complex
Function Theory}, Springer, Berlin, 1998.

\item[{\bf [SIC]}] J. Siciak, On series of homogeneous
polynomials and their partial sums, {\it Annales Polonici
Math.} {\bf 51}(1990), 289--302.

\item[{\bf [STE]}] L. Stepnickova, Pointwise and locally
uniform convergence of holomorphic and harmonic functions,
{\it Comment.\ Math.\ Univ.\ Carolin.} 40(1999), 665--678.

\item[{\bf [STE]}] T.-J. Stieltjes, Recherches sur les
fractions continues, {\it Ann.\ Fac.\ Sci.\ Toulouse}
VIIIJ(1894), 1--122.

\item[{\bf [TIT]}]  E. C. Titchmarsh, {\it The Theory of Functions}, Oxford
University Press, New York, 1960.

\item[{\bf [WEI]}]  K. Weierstrass, {\it Mathematische Werke}, Hildesheim, G. Olms,
New York, Johnson Reprint, 1967.

\item[{\bf [ZAL]}]  L. Zalcman, Real proofs of complex theorems (and vice versa),
{\it Am.\ Math.\ Monthly} 81(1974), 115--137.

\end{enumerate}
\vspace*{.25in} 

\noindent Department of Mathematics, Washington University in St.\ Louis, St. Louis, Missouri, USA 63130
\smallskip  \\
{\tt sk@math.wustl.edu}

\end{document}